\documentclass[11pt]{amsart}
\usepackage{amscd, amssymb} 
\input{prepictex.tex}
\input{pictex.tex}
\input{postpictex.tex}

\def\T{\mathbb T}

\def\C{\mathbb C}

\def\H{{\mathcal H}}
\def\K{{\mathcal K}}

\def\O{{\mathcal O}}

\def\Kbad{K^{\text{{\rm fin}}}_{\infty}}

\def\KMbad{(K _M)^{\text{{\rm fin}}}_{\infty}}
\def\Kppbad{(K'')^{\text{{\rm fin}}}_{\infty}}
\def\Kpbad{(K')^{\text{{\rm fin}}}_{\infty}}

\def\OMbad{\Omega(M)^{\text{{\rm fin}}}_{\infty}}

\def\ONbad{\Omega(N)^{\text{{\rm fin}}}_{\infty}}

\def\Ovbad{\Omega(v)^{\text{{\rm fin}}}_{\infty}}
\def\OUbad{\Omega(U)^{\text{{\rm fin}}}_{\infty}}
\def\OUpbad{\Omega(U')^{\text{{\rm fin}}}_{\infty}}
\def\Owbad{\Omega(w)^{\text{{\rm fin}}}_{\infty}} 

\def\M0{\Omega(M)_\infty^\emptyset} 

\def\Aut{\operatorname{Aut}}

\def\id{\text{id}}

\def\Prim{\operatorname{Prim}}

\newtheorem{theo}{Theorem}[section]
\newtheorem{cor}[theo]{Corollary}
\newtheorem{lemma}[theo]{Lemma}

\theoremstyle{definition}
\newtheorem{definition}[theo]{Definition}

\theoremstyle{remark}
\newtheorem{remark}[theo]{Remark}
\newtheorem{example}[theo]{Example}

\begin{document}
\title[Ideals of Graph Algebras]{The Primitive Ideal Space 
       of the $C^*$-Algebras of Infinite Graphs} 

\author[Hong]{Jeong Hee Hong}
\address{Applied Mathematics\\
Korea Maritime University\\
Pusan 606--791\\
South Korea}
\email{hongjh@hanara.kmaritime.ac.kr}

\author[Szyma\'{n}ski]{Wojciech Szyma\'{n}ski}
\address{Mathematics\\
The University of Newcastle\\
Callaghan, NSW 2308\\
Australia}
\email{wojciech@frey.newcastle.edu.au}

\date{} 

\maketitle

\begin{abstract}
For any countable directed graph $E$ we describe the primitive ideal 
space of the corresponding generalized Cuntz-Krieger algebra $C^*(E)$. 
\end{abstract}

\addtocounter{section}{-1} 

\section{Introduction} 

The primary purpose of this article is to give a description of 
the primitive ideal space of the $C^*$-algebra $C^*(E)$ corresponding 
to an arbitrary countable directed graph $E$. The two main results of 
the article, Theorem \ref{otherprimitive} (together with Corollary 
\ref{primitiveideals}) and Theorem \ref{9closures}, identify elements 
of $\Prim(C^*(E))$ and describe the closure operation in the hull-kernel 
topology. These theorems build on and generalize a long string of previous 
results on the ideal structure of Cuntz-Krieger type algebras, obtained 
by a number of researchers over the period of past twenty years. Our 
present article completes the program of classification of ideals 
of the generalized Cuntz-Krieger algebras corresponding to arbitrary 
countable directed graphs. 

First fundamental results about the ideal structure of Cuntz-Krieger 
algebras were obtained by Cuntz, who described all ideals 
of $\O_A$ for a finite $0-1$ matrix $A$ satisfying Condition (II) 
\cite[Theorem 2.5]{c}. Much more recently, an Huef and Raeburn gave 
a complete description of all gauge-inaviant ideals \cite[Theorem 3.5]{hr} 
and the primitive ideal space \cite[Theorem 4.7]{hr} for Cuntz-Krieger 
algebras $\O_A$ corresponding to arbitrary finite $0-1$ matrices. 
Soon after the Cuntz-Krieger algebras of countably infinite directed 
graphs (the graph algebras) were introduced and analyzed by Kumjian, 
Pask, Raeburn and Renault. They described all ideals of $C^*(E)$ for a 
locally finite graph $E$ satisfying Condition (K) (an analogue of 
Cuntz's Condition (II)) \cite[Theorem 6.6]{kprr}. Since then a number 
of papers considered the problem of classification of ideals of 
graph algebras and other generalizations of the classical Cuntz-Krieger 
algebras. However, to the best of our knowledge, most of those papers 
dealt only with ideals invariant under the gauge action. This in 
particular applies to graphs satisfying Condition (K), since 
for such graphs all closed ideals of $C^*(E)$ are gauge-invariant 
\cite[Lemma 2.2]{hs3}. For row-finite graphs, Bates, Pask, Raeburn 
and Szyma\'nski gave all gauge-invariant ideals of $C^*(E)$ 
\cite[Theorem 4.1]{bprs} and described the primitive ideal space if 
in addition $E$ satisfies Condition (K) \cite[Theorem 6.3]{bprs}. 
Working with arbitrary countable graphs, Bates, Hong, Raeburn and 
Szyma\'nski described all gauge-invariant ideals \cite[Theorem 3.6]{bhrs} 
and identify those of them which are primitive \cite[Theorem 4.7]{bhrs}. 
(A brief overview of the results of \cite{bhrs} was reported by Hong 
in \cite{h}.) About the same time, Drinen and Tomforde obtained 
similar results (through different techniques) for graphs satisfying 
Condition (K) \cite[Theorem 3.5 and Theorem 4.10]{dt}. 

Besides graph algebras, there are other interesting generalizations 
of Cuntz-Krieger algebras. These include Cuntz-Pimsner algebras 
generated by Hilbert bimodules and Exel-Laca algebras  
corresponding to infinite $0-1$ matrices. Many, though not all, 
graph algebras can be viewed as either Cuntz-Pimsner algebras 
or Exel-Laca algebras (cf. \cite{fr,sz0}). Partial 
results about the ideal structure of Cuntz-Pimsner algebras 
were obtained by Pinzari \cite{pin}, 
Kajiwara, Pinzari and Watatani \cite{kpw}, and Fowler, Muhly and 
Raeburn \cite{fmr}. Exel and Laca described ideals of the Cuntz-Krieger 
algebras $\O_A$ corresponding to infinite $0-1$ matrices, under 
an extra hypothesis on the matrix $A$ analogous to Condition (K) 
\cite[Theorem 15.1]{el}. 

Obviously, there are many benefits from such a comprehensive 
description of the ideal structure of a large class of algebras, 
as presented in \cite{bhrs} and the present article. For example, 
Szyma\'nski's proof of a very general criterion of injectivity of 
homomorphisms of graph algebras \cite[Theorem 1.2]{sz2} is based 
on the analysis of ideals, and so are some arguments from the recent 
work of Hong and Szyma\'nski on non-simple purely infinite graph algebras 
\cite{hs3}. In fact, we think that graph algebras might play 
a prominent role in the study and classification of this interesting 
class of $C^*$-algebras, whose investigations have been recently 
initiated by Kirchberg and R\o rdam \cite{kr1,kr2}. 
Also, graph algebras appearing in the context of some 
compact quantum manifolds \cite{hs1,hs2} are not simple, and it is 
important to know their ideal structure. Furthermore, it does not 
seem unlikely, that our methods, techniques and reults 
on the ideal structure of graph algebras may help in understanding of 
other classes of $C^*$-algebras, Cuntz-Pimsner algebras for example. 
Indeed, techniques quite similar to those developed 
for the study of graph algebras have recently been used by 
Katsura in his analysis of the crossed products of Cuntz algebras 
by quasi-free actions of locally compact abelian groups \cite{k1,k2}. 
And certainly, good understanding of the ideal structure of 
generalized Cuntz-Krieger algebras is a necessary first step 
towards their classification. 

Of course, as a by-product of our analysis of the ideal 
structure of graph algebras we obtain their simplicity criteria.  
The problem of simplicity of generalized Cuntz-Krieger algebras 
was discussed by a number of authors. Partial answers to this 
question for graph algebras were given by Bates, Pask, Raeburn 
and Szyma\'nski in \cite[Proposition 5.1]{bprs}, and by Fowler,  
Laca and Raeburn in \cite[Theorem 3]{flr}. An if and only if 
criterion was proved by Szyma\'nski in \cite[Theorem 12]{sz1}, 
and another one somewhat later but independently by Paterson in 
\cite[Theorem 4]{pat}. A partial result about simplicity of 
$\O_A$ for an infinite matrix $A$ was given by Exel and Laca in 
\cite[Theorem 14.1]{el}, and an if and only if criterion was 
supplied by Szyma\'nski in \cite[Theorem 8]{sz1}. Building on 
this latter result, Tomforde proved an analogous criterion for 
the $C^*$-algebras corresponding to ultragraphs \cite{t} 
(a class of $C^*$-algebras which contains both graph algebras 
$C^*(E)$ and Exel-Laca algebras $\O_A$). Related results about 
simplicity of Cuntz-Pimsner algebras were obtained by Schweizer 
in \cite{sch1} and \cite{sch2}. 

The present article is organized as follows. In \S1, we review basic 
facts about graph algebras we need. The main reference to this section 
is \cite{bhrs}. We rely heavily on the results of that paper. In particular, 
we use the description of quotients of $C^*(E)$ by gauge-invariant 
ideals \cite[Proposition 3.4]{bhrs}, the description of the 
intersection of a family of gauge-invariant ideals 
\cite[Proposition 3.9 and Corollary 3.10]{bhrs}, 
the concepts of maximal tails and breaking vertices \cite[\S4]{bhrs}, 
the description of all gauge-invariant ideals of 
$C^*(E)$ \cite[Theorem 3.6]{bhrs}, 
and the classification of gauge-invariant primitive ideals of $C^*(E)$ 
\cite[Theorem 4.7]{bhrs}. In \S2, we describe all primitive ideals of 
$C^*(E)$ which are not invariant under the gauge action (cf. Theorem 
\ref{otherprimitive}). Inside $\Prim(C^*(E))$ they form circles which are 
in one-to-one correspondence with maximal tails containing  
a loop without exits (cf. Lemma \ref{looptail}). The general plan of our 
argument is similar to that of \cite{hr} and is based on sandwiching a 
non gauge-invariant primitive ideal between two gauge-invariant ideals 
(cf. Lemmas \ref{sandwich} and \ref{loopquotient}). However, the case 
of arbitrary infinite graphs is technically much more complicated than 
that of finite graphs. The main result of \S2 is Corollary 
\ref{primitiveideals}, which gives a description of all primitive 
ideals of $C^*(E)$ for an arbitrary countable graph $E$. In \S3, we 
describe the closure operation in the hull-kernel topology of 
$\Prim(C^*(E))$ (cf. Theorem \ref{9closures}). Since our result covers 
all possible countable directed graphs, this description is necessarily 
somewhat involved. It greatly simplifies in the case of row-finite 
graphs (cf. Corollary \ref{4closures}). In \S4, we illustrate the 
main results with a few examples. 

\vspace{2mm}\noindent 
{\bf Acknowledgements.} We are grateful to Professor Iain Raeburn for 
many a valuable discussion. His much helpful advice stimulated and 
influenced our entire line of investigation. The first named author 
would like to thank the faculty and staff of the Mathematics Department 
at the Newcastle University for their warm hospitality during her several 
visits there. The second named author wishes to thank members of the 
Mathematics Department at the Korea Maritime University for their 
hospitality during his stay there. He is  grateful to Professor Hideki 
Kosaki and the Japanese Society for the Promotion of Science for their 
generous support. He is indebted to the Mathematics Departments of Kyushu 
and Kyoto Universities, and especially to Professor Masaki Izumi, for 
their kind hospitality. 

\section{Preliminaries on graph algebras}

We recall the definition of the $C^*$-algebra 
corresponding to a directed graph \cite{flr}. 
Let $E=(E^0,E^1,r,s)$ be a directed graph with 
countably many vertices $E^0$ and  
edges $E^1$, and range and source functions 
$r,s:E^1\rightarrow E^0$, respectively. 
$C^*(E)$ is defined as the universal 
$C^*$-algebra generated by families of projections 
$\{p_v:v\in E^0\}$ and partial isometries $\{s_e:e\in E^1\}$, 
subject to the following relations.  
\begin{description} 
\item[(GA1)] $p_vp_w=0$ for $v,w\in E^0$, $v\neq w$. 
\item[(GA2)] $s_e^*s_f=0$ for $e,f\in E^1$, $e\neq f$. 
\item[(GA3)] $s_e^*s_e=p_{r(e)}$ for $e\in E^1$.  
\item[(GA4)] $s_e s_e^*\leq p_{s(e)}$ for $e\in E^1$.  
\item[(GA5)] $p_v=\displaystyle{\sum_{e\in E^1:\;
s(e)=v}}s_e s_e^*$ for $v\in E^0$ such that $0<|s^{-1}(v)|<\infty$. 
\end{description} 
Universality in this definition means that if $\{Q_v:v\in E^0\}$ and 
$\{T_e:e\in E^1\}$ are families of projections and partial 
isometries, respectively, satisfying conditions (GA1--GA5), then 
there exists a $C^*$-algebra homomorphism from $C^*(E)$ to the 
$C^*$-algebra generated by $\{Q_v:v\in E^0\}$ and 
$\{T_e:e\in E^1\}$ such that $p_v\mapsto Q_v$ and 
$s_e\mapsto T_e$ for $v\in E^0$, $e\in E^1$. 

It follows from the universal property that there exists a gauge 
action $\gamma:\T\rightarrow\Aut(C^*(E))$ such that 
$\gamma_z(p_v)=p_v$ and $\gamma_z(s_e)=zs_e$ for all 
$v\in E^0$, $e\in E^1$, $z\in\T$. We denote by $\Gamma$ the 
corresponding conditional expectation of $C^*(E)$ onto the 
fixed-point algebra $C^*(E)^\gamma$, such that 
$\Gamma(x)=\int_{z\in\T}\gamma_z(x)dz$ for $x\in C^*(E)$. The integral 
is with respect to the normalized Haar measure on $\T$. Note that 
$\Gamma(p_v)=p_v$ and $\Gamma(s_e)=0$ for all $v\in E^0$, $e\in E^1$. 

If $\alpha_1,\ldots,\alpha_n$ are  (not necessarily  distinct) 
edges such that $r(\alpha_i)= s(\alpha_{i+1})$ for $i=1,\ldots,n-1$, 
then $\alpha=(\alpha_1,\ldots,\alpha_n)$ is a path of length 
$|\alpha|=n$, with source $s(\alpha)=s(\alpha_1)$ and range  
$r(\alpha)=r(\alpha_n)$. We also allow paths of length zero, 
identified with vertices. The set of all paths of length $n$ 
is denoted by $E^n$, while the collection of all finite paths in 
$E$ is denoted by $E^*$. Given a path $\alpha=(\alpha_1,\ldots,\alpha_n)$ 
we denote $s_\alpha=s_{\alpha_1}\cdots s_{\alpha_n}$, a partial 
isometry in $C^*(E)$. A loop is a path of positive 
length whose source and range coincide. A loop $\alpha$
has an exit if there exists an edge $e\in E^1$ and index $i$
such that $s(e)=s(\alpha_i)$ but $e\neq\alpha_i$. If $\alpha$ is a loop 
all of whose vertices belong to a subset $M\subseteq E^0$ then we say 
that $\alpha$ has an exit in $M$ if an edge $e$ exists as above with 
$r(e)\in M$. 

By an ideal in a $C^*$-algebra we always mean a closed two-sided ideal. 
In order to understand the ideal structure of a graph algebra it is 
convenient to look at saturated hereditary subsets of the vertex set. 
As usual, if $v,w\in E^0$ then we write $v\geq w$ when
there is a path from $v$ to $w$, and say that a subset $K$ of $E^0$ is 
hereditary if $v\in K$ and $v\geq w$ imply $w\in K$.  A subset
$K$ of $E^0$ is saturated if every vertex $v$ which  satisfies
$0<|s^{-1}(v)|<\infty$ and $s(e)=v\Longrightarrow r(e)\in K$ 
itself belongs to $K$. If $X\subseteq E^0$ then $\Sigma(X)$ is the 
smallest saturated subset of $E^0$ containing $X$, and $\Sigma H(X)$ 
is the smallest saturated hereditary subset of $E^0$ containing $X$. 
If $K$ is hereditary and saturated then $I_K$ denotes the ideal 
of $C^*(E)$ generated by $\{p_v:v\in K\}$. We have 
$$ I_K=\overline{\text{span}}\{s_\alpha p_v s_\beta^*:
   \alpha,\beta\in E^*,\:v\in K,\:r(\alpha)=r(\beta)=v\}. $$ 
As shown in \cite[Proposition 3.4]{bhrs}, the quotient $C^*(E)/I_K$ 
is naturally isomorphic to $C^*(E/K)$. The quotient graph $E/K$ was  
defined in \cite[Section 3]{bhrs}. The vertices of $E/K$ are 
$$ (E/K)^0=(E^0\setminus K)\cup\{\beta(v):v\in \Kbad\}, $$ 
where 
$$ \Kbad=\{v\in E^0\setminus K:|s^{-1}(v)|=\infty\text{ and } 
   0<|s^{-1}(v)\cap r^{-1}(E^0\setminus K)|<\infty\}. $$ 
The edges of $E/K$ are 
$$ (E/K)^1=r^{-1}(E^0\setminus K)\cup\{\beta(e):
   e\in E^1,r(e)\in\Kbad\}, $$ 
with $r,s$ extended by $r(\beta(e))=\beta(r(e))$ and 
$s(\beta(e))=s(e)$. Note that all extra vertices $\beta(\Kbad)$ are sinks 
in $E/K$. If $\Kbad=\emptyset$ then $E/K$ is simply a subgraph of $E$, 
denoted $E\setminus K$. If $v\in\Kbad$ then we write 
$$ p_{v,K}=\sum_{s(e)=v,\,r(e)\not\in K}s_es_e^*. $$  
For $B\subseteq\Kbad$ the ideal of $C^*(E)$ 
generated by $I_K$ and $\{p_v-p_{v,K}:v\in B\}$ is denoted by $J_{K,B}$. 
We have 
$$ J_{K,B}=\overline{\text{span}}\big\{s_\alpha p_v s_\beta^*,\;
   s_\mu(p_w-p_{w,K})s_\nu^*:\alpha,\beta,\mu,\nu\in E^*, $$ 
$$ \;\;\;\;\;\;\;\;\;\;\;\;\;\;\;\;\;\;\;\;\;\;\;\;\;\;\;\;\;\;\;\;\;
   v\in K,w\in B,r(\alpha)=r(\beta)=v,\: r(\mu)=r(\nu)=w\big\}. $$ 
By \cite[Corollary 3.5]{bhrs}, the quotient $C^*(E)/J_{K,B}$ is 
naturally isomorphic to $C^*((E/K)\setminus\beta(B))$. 
As shown in \cite[Theorem 3.6]{bhrs}, all gauge-invariant ideals of 
$C^*(E)$ are of the form $J_{K,B}$. 

A non-empty subset $M\subseteq E^0$ is a maximal tail if it satisfies 
the following three conditions (cf. \cite[Lemma 4.1]{bhrs}): 
\begin{description} 
\item{(MT1)} If $v\in E^0$, $w\in M$, and $v\geq w$, then $v\in M$.  
\item{(MT2)} If $v\in M$ and $0<|s^{-1}(v)|<\infty$, then there 
exists $e\in E^1$ with $s(e)=v$ and $r(e)\in M$. 
\item{(MT3)} For every $v,w\in M$ there exists $y\in M$ such that 
$v\geq y$ and $w\geq y$. 
\end{description}  
We denote by ${\mathcal M}(E)$ the collection of all maximal tails 
in $E$ and by ${\mathcal M}_\gamma(E)$ the collection of all maximal 
tails $M$ in $E$ such that each loop in $M$ has an exit in $M$. We set 
${\mathcal M}_\tau(E)={\mathcal M}(E)\setminus{\mathcal M}_\gamma(E)$.   

If $X\subseteq E^0$ then, as in \cite{bhrs}, we denote  
$$ \Omega(X)=\{w\in E^0\setminus X:w\not\geq v 
   \text{ for all } v\in X\}. $$ 
If $X$ consists of a single vertex $\{v\}$ then we write 
$\Omega(v)$ instead of $\Omega(\{v\})$. 
Note that $\Omega(M)=E^0\setminus M$ for every maximal tail $M$. Moreover, 
for such an $M$, $\Omega(M)$ is hereditary by (MT1) and saturated by (MT2). 

Along with maximal tails, the set 
$$ BV(E)=\{v\in E^0:|s^{-1}(v)|=\infty \text{ and } 0<|s^{-1}(v)
   \setminus r^{-1}(\Omega(v))|<\infty\} $$ 
plays an important role in the classification of primitive gauge-invariant 
ideals \cite{bhrs}. Its elements are called breaking vertices. Note that 
if $K\subseteq E^0$ is hereditary and saturated, then $v\in\Kbad$ implies  
$v\in BV(E)$. If $v\in BV(E)$ then $\Omega(v)$ is hereditary and saturated. 

We denote by $\Prim(C^*(E))$ the set of all primitive ideals of $C^*(E)$  
and by $\Prim_\gamma(C^*(E))$ the set of all primitive gauge-invariant 
ideals of $C^*(E)$.  

As shown in \cite[Theorem 4.7]{bhrs}, there is a one-to-one correspondence  
$$ {\mathcal M}_\gamma(E)\cup BV(E)\longrightarrow\Prim_\gamma(C^*(E)) $$
given by 
\begin{align*} 
{\mathcal M}_\gamma(E)\ni M & \longmapsto  J_{\Omega(M),\OMbad}, \\ 
BV(E)\ni v & \longmapsto  J_{\Omega(v),\Ovbad\setminus\{v\}}. 
\end{align*} 

\section{The primitive ideals} 

Our goal in this section is to show that any maximal tail 
in ${\mathcal M}_\tau(E)$ gives rise to 
a circle of primitive ideals, none of which is gauge-invariant, 
and that all non gauge-invariant primitive 
ideals arise in this way. To this 
end we explicitly construct the corresponding irreducible 
representations of $C^*(E)$. At first we observe that any maximal 
tail  $M\in{\mathcal M}_\tau(E)$ contains an essentially unique 
vertex simple loop without exits in $M$. A loop $L=(e_1,\ldots,e_n)$ 
is vertex simple if and only if $r(e_i)\neq r(e_j)$ for $i\neq j$. 
We denote by $L^0$ the set $\{r(e_i):i=1,\ldots,n\}$ of the vertices 
through which $L$ passes, and by $L^1$ the set $\{e_i:i=1,\ldots,n\}$ 
of its edges. 

\begin{lemma}\label{looptail}
If $M$ is a maximal tail in $E$ then $M\in{\mathcal M}_\tau(E)$ if and 
only if there exists a vertex simple loop $L$ with $L^0\subseteq M$ and 
such that if $e\in E^1\setminus L^1$ and $s(e)\in L^0$ then $r(e)\not\in M$.  
Furthermore, such a loop is unique up to a cyclic permutation 
of the edges comprising it, and $\Omega(M)=\Omega(L^0)$. 
\end{lemma} 
\begin{proof} 
The first assertion follows immediately from the 
definition of ${\mathcal M}_\tau(E)$ as the complement 
of ${\mathcal M}_\gamma(E)$ in ${\mathcal M}(E)$. 

For the uniqueness suppose that there are 
two loops $L_1=(e_1,\ldots,e_n)$ and $L_2=(f_1,\ldots,f_m)$ 
with the above properties. By condition (MT3) there 
is a vertex $v\in M$ and paths $\alpha,\beta$ such that $s(\alpha)=s(e_1)$, 
$s(\beta)=s(f_1)$, and $r(\alpha)=r(\beta)=v$. Since $L_1$, $L_2$ have no 
exits in $M$ we must have $v=s(e_k)=s(f_r)$ for some $k,r$. The absence 
of exits then implies that $e_k=f_r$, $e_{k+1}=f_{r+1}$, and so on, which 
proves the claim. 

Obviously $\Omega(M)\subseteq\Omega(L^0)$. For the reverse inclusion 
let $L=(e_1,\ldots,e_n)$ and $v\in\Omega(L^0)$. Suppose for a moment that 
$v\in M$. By (MT3), there exists a vertex $w\in M$ such that 
$v\leq w$ and $s(e_1)\leq w$. Since $L$ has no exits in $M$, $w$ 
must be in $L^0$. Thus $v\not\in\Omega(L^0)$, contrary to the assumption. 
Hence $v\not\in M$. Since $M$ is a maximal tail we have $E^0\setminus M=
\Omega(M)$ and hence $v\in\Omega(M)$, as required. 
\end{proof} 

From now on, for each maximal tail 
$M\in{\mathcal M}_\tau(E)$ we choose one vertex simple loop 
without exits in $M$ as in Lemma \ref{looptail} and call it $L_M$. 
For $v\in E^0\setminus L_M^0$ we denote 
$$ A_M(v)=\{(\alpha_1,\ldots,\alpha_m)\in E^*:s(\alpha_1)=v,\:r(\alpha_m)
   \in L_M^0,\:r(\alpha_i)\not\in L_M^0 \text{ if } i\neq m\} $$ 
and set 
$$ A_M=\bigcup_{v\in E^0\setminus L_M^0}A_M(v). $$ 

\begin{definition}\label{L_Mloop} 
Let $E$ be a directed graph. Given $M\in{\mathcal M}_\tau(E)$ let  
$L_M=(e_1,\ldots,e_n)$. We denote 
by $\H_M$ the Hilbert space with an orthonormal basis $\{\xi_\alpha
:\alpha\in A_M\cup L_M^0\}$. For $t\in\T\subset\C$ we define a 
representation 
$$ \rho_{M,t}:C^*(E)\longrightarrow{\mathcal B}(\H_M) $$ 
so that  
\begin{align*}
\rho_{M,t}(p_v)\xi_\alpha & =  \left\{ \begin{array}{ll} 
     \xi_\alpha & \mbox{ if } s(\alpha)=v \\ 
     0 & \mbox{ otherwise, } \end{array} \right. \\ 
\rho_{M,t}(s_e)\xi_\alpha & =  \left\{ \begin{array}{ll} 
     t\xi_{s(e)} & \mbox{ if } e=e_1 \mbox{ and } \alpha=r(e_1) \\ 
     \xi_{s(e)} & \mbox{ if } r(e)=\alpha\in L_M^0 \mbox{ and } e_1\neq e\in L_M^1 \\ 
     \xi_e & \mbox{ if } r(e)=\alpha\in L_M^0 \mbox{ and } e\not\in L_M^1 \\ 
     \xi_{(e,\alpha_1,\ldots,\alpha_m)} & \mbox{ if } 
     \alpha,(e,\alpha_1,\ldots\alpha_m)\in A_M \\ 
     0 & \mbox{ otherwise, } \end{array}\right. 
\end{align*} 
for $v\in E^0$, $e\in E^1$, and $\alpha=
(\alpha_1,\ldots,\alpha_m)\in A_M\cup L_M^0$. 
\end{definition} 

\begin{remark}\label{abuse}  
That $\rho_{M,t}$ indeed gives rise to a representation of 
$C^*(E)$ will be shown in Lemma \ref{rhoirreducible}. 
Strictly speaking, this representation depends not only on the 
maximal tail $M$ but also on the choice of the loop $L_M$. We 
slightly abuse the notation by writing $\rho_{M,t}$ instead of 
more precise $\rho_{L_M,t}$. We will use the latter notation later 
in Lemma \ref{loopquotient} for emphasis, when considering $L_M$ and 
its cyclic permutation $L'_M$ simultaneously. 
\end{remark} 

We will see (cf. Lemma \ref{rhoirreducible}) that each 
$\rho_{M,t}$ is an irreducible representation of $C^*(E)$ 
and thus $\ker\rho_{M,t}$ is a primitive ideal, which turns out 
to be not invariant under the gauge action. It will be useful to 
sandwich  such an ideal between two gauge-invariant ones (cf. Lemma 
\ref{sandwich}), and to this end we need to consider the set 
$$ K_M:=L_M^0\cup\{v\in E^0\setminus L_M^0:|A_M(v)|<\infty\}. $$  
It is not difficult to see that each $K_M$ is 
hereditary, saturated and that $\Sigma H(\Omega(M)\cup L_M^0)
\subseteq K_M$. If $E$ is row-finite then $\Sigma H(\Omega(M)\cup L_M^0)
=K_M$ but in general $K_M$ may be larger. 

\begin{lemma}\label{projectionsinKM} 
Let $E$ be a directed graph, $M\in{\mathcal M}_\tau(E)$ and   
$\pi:C^*(E)\rightarrow C^*(E)/J_{\Omega(M),\OMbad}$ be  
the natural surjection. Then for every vertex $v\in K_M\setminus
(\Omega(M)\cup L_M^0)$ we have  
$$ \pi(p_v)=\sum_{\alpha\in A_M(v)}\pi(s_\alpha s_{\alpha}^*). $$ 
\end{lemma} 
\begin{proof} 
By \cite[Corollary 3.5]{bhrs}, there exists a natural isomorphism 
between $C^*(E)/J_{\Omega(M),\OMbad}$ and $C^*(E\setminus\Omega(M))$.  
Since $E\setminus\Omega(M)=(M,r^{-1}(M),r,s)$  
we may assume that $M=E^0$ and $\pi=\id$. Now 
we proceed by induction with respect to the maximum length $\ell$ 
of elements of $A_M(v)$. If $\ell=1$ then the claim follows 
from condition (GA5). For the inductive step observe that  
if $v\in K_M\setminus(\Omega(M)\cup L_M^0)$ then the set 
$\{e\in E^1:s(e)=v\}$ is non-empty by Lemma~\ref{looptail} 
and finite by the definition of $K_M$. Thus 
$p_v=\sum_{e\in E^1,s(e)=v}s_es_e^*$ by (GA5). 
Applying the inductive hypothesis to 
$\{p_w:w=r(e),\:e\in E^1,\:s(e)=v\}$ 
we infer that the desired identity holds. 
\end{proof} 

If the graph $E$ is row-finite then the smallest gauge-invariant 
ideal of $C^*(E)$ containing $\ker\rho_{M,t}$ is $I_{K_M}$ 
(cf. Lemma \ref{sandwich}). However, if $E$ contains vertices with 
infinite valencies then such a gauge-invariant ideal must have 
the form $J_{K_M,B_M}$ for a suitable $B_M\subseteq\KMbad$. 
It turns out that  
$$ B_M:=\KMbad\cap\OMbad $$ 
does the trick. Note that a vertex $v\in E^0$ belongs to 
$B_M$ if and only if among $\{e\in E^1:s(e)=v\}$ there are 
infinitely many edges $e$ such that $r(e)\in\Omega(M)$, only finitely 
many $e$ with $r(e)\in E^0\setminus\Omega(M)=M$, and at 
least one $e$ such that $r(e)\not\in K_M$. 

\begin{lemma}\label{rhoirreducible} 
Let $E$ be a directed graph, $M\in{\mathcal M}_\tau(E)$ 
and $t\in\T$. Then $\rho_{M,t}$ of Definition \ref{L_Mloop} 
gives rise to an irreducible representation 
of $C^*(E)$ such that $\rho_{M,t}(J_{K_M,B_M})=\K(\H_M)$. 
\end{lemma} 
\begin{proof} At first we show that $\{\rho_{M,t}
(p_v),\rho_{M,t}(s_e):v\in E^0,e\in E^1\}$ is a Cuntz-Krieger 
$E$-family. Obviously, conditions (GA1--GA4) are satisfied. We 
verify (GA5). Let $v$ be a vertex such that $0<|s^{-1}(v)|
<\infty$. If $v\in L_M^0$ then there is exactly 
one edge $f$ with $s(f)=v$ and $r(f)\in L_M^0$. It then 
follows from the definition of $\rho_{M,t}$ that  
both $\sum_{e\in E^1,s(e)=v}\rho_{M,t}(s_e)\rho_{M,t}(s_e)^*=\rho_{M,t}
(s_f)\rho_{M,t}(s_f)^*$ and $\rho_{M,t}(p_v)$ are projections onto the 
1-dimensional subspace of $\H_M$ spanned by $\xi_v$. If $v\in E^0\setminus 
L_M^0$ then $\rho_{M,t}(p_v)$ is a projection onto the subspace of $\H_M$ 
spanned by $\{\xi_\alpha:\alpha\in A_M(v)\}$. If  
$e\in E^1$, $s(e)=v$, then $\rho_{M,t}(s_e)^*\xi_\alpha=0$ for  
$s(\alpha)\neq v$. If $s(\alpha)=v$ and $\alpha=(f_1,\ldots,f_k)$, 
$f_i\in E^1$, then $\rho_{M,t}(s_e)^*\xi_\alpha=0$ for $e\neq f_1$ and 
$\rho_{M,t}(s_e)\rho_{M,t}(s_e)^*\xi_\alpha=
\xi_\alpha$ for $e=f_1$. Thus again 
$\sum_{e\in E^1,s(e)=v}\rho_{M,t}
(s_e)\rho_{M,t}(s_e)^*=\rho_{M,t}(p_v)$.  
Hence, Definition \ref{L_Mloop} gives rise to a 
representation $\rho_{M,t}:C^*(E)\rightarrow{\mathcal B}(\H_M)$. 

It follows from Definition \ref{L_Mloop} that all the projections 
$\rho_{M,t}(p_v)$, $v\in K_M$, and $\rho_{M,t}(p_w-p_{w,K_M})$, 
$w\in B_M$, are of finite rank. Thus $\rho_{M,t}(J_{K_M,B_M})
\subseteq\K(\H_M)$. On the other hand, the range $\rho_{M,t}(C^*(E))$ 
contains $\K(\H_M)$. Indeed, if $v\in L_M^0$ then $\{\rho_{M,t}
(s_\mu):\mu\in E^*,\:r(\mu)=v\}$ are rank-1 partial 
isometries sending $\xi_v$ to all other elements of the 
orthonormal basis $\{\xi_\alpha:\alpha\in A_M\cup L_M^0\}$. 
Consequently $\rho_{M,t}(J_{K_M,B_M})=\K(\H_M)$. In particular, 
$\rho_{M,t}$ is irreducible. 
\end{proof} 

By Lemma \ref{rhoirreducible}, the representations $\rho_{M,t}$ 
give rise to primitive ideals $\ker\rho_{M,t}$ 
of $C^*(E)$. In Lemmas \ref{sandwich} and \ref{loopquotient}, 
below, we show the key property of $C^*(E)$ 
that for each $M\in{\mathcal M}_\tau(E)$ the union of 
$\{\ker\,\rho_{M,t}:t\in\T\}$ may be sandwiched between two 
uniquely determined gauge-invariant ideals whose quotient is 
Morita equivalent to $C(\T)$. This result was originally 
proved by an Huef and Raeburn in \cite[Lemma 4.5]{hr} for the 
Cuntz-Krieger algebras corresponding to finite matrices. The 
argument there took advantage of the existence   
of only finitely many gauge-invariant ideals. In the present 
article we need a different argument, as algebras corresponding to 
infinite graphs may have infinitely many gauge-invariant ideals. 

\begin{lemma}\label{sandwich} 
Let $E$ be a directed graph, $M\in{\mathcal M}_\tau(E)$  
and $t\in\T$. Then the following hold. 
\begin{enumerate} 
\item The ideal $J_{\Omega(M),\OMbad}$ is the largest among 
gauge-invariant ideals of $C^*(E)$ contained in $\ker\rho_{M,t}$. 
\item The ideal $J_{K_M,B_M}$ is the smallest among 
gauge-invariant ideals of $C^*(E)$ containing $\ker\rho_{M,t}$.     
\end{enumerate} 
\end{lemma} 
\begin{proof} 
Ad 1. 
Since $\Omega(M)=\Omega(L_M^0)$ by Lemma \ref{looptail}, 
it is immediate from Definition \ref{L_Mloop} that 
$\rho_{M,t}(p_v)=0$ if $v\in\Omega(M)$, and $\rho_{M,t} 
(p_v-p_{v,\Omega(M)})=0$ if $v\in \OMbad$. Hence 
$J_{\Omega(M),\OMbad}\subseteq\ker\rho_{M,t}$. 

Let $J_1$ be a gauge-invariant ideal of $C^*(E)$ 
contained in $\ker\rho_{M,t}$. By \cite[Theorem 3.6]{bhrs}  
there is a saturated hereditary $K\subseteq E^0$ and a 
$B\subseteq \Kbad$ such that $J_1=J_{K,B}$. By 
\cite[Corollary 3.10]{bhrs}, in order that 
$J_{K,B}\subseteq J_{\Omega(M),\OMbad}$ we must 
have $K\subseteq\Omega(M)$ and $B\subseteq \Omega(M)\cup 
\OMbad$. Let $v\in E^0\setminus\Omega(M)$. 
Since $\Omega(M)=\Omega(L_M^0)$ there 
is a path from $v$ to $L_M^0$. If $\alpha$ is such a path with 
the shortest possible length then $\alpha\in A_M\cup L_M^0$ and 
$\rho_{M,t}(p_v)\xi_\alpha\neq0$. Thus $\rho_{M,t}(p_v)\neq0$ and  
consequently $v\not\in K$. This shows that $K\subseteq\Omega(M)$. 
Now let $v\in B\setminus\Omega(M)$. Then $v$ emits at least one edge 
into $M$, since there is a path from $v$ to $L_M^0$. Furthermore, $v$ emits 
infinitely many edges into $K\subseteq\Omega(M)$, and only finitely 
many edges into $E^0\setminus K$ and hence into $E^0\setminus\Omega(M)$. 
Consequently $v\in \OMbad$. Thus $B\subseteq\Omega(M)\cup\OMbad$. 

Ad 2. Let $J=J_{K_M,B_M}$. 
We denote by $\pi:{\mathcal B}(\H_M)\rightarrow{\mathcal B}
(\H_M)/\K(\H_M)$ the quotient map of ${\mathcal B}(\H_M)$ onto 
its Calkin algebra. By Lemma \ref{rhoirreducible} we have 
$\rho_{M,t}(J)=\K(\H_M)$. 

At first we show that $\ker\rho_{M,t}\subseteq J$. 
To this end it suffices to prove injectivity of the homomorphism 
$$ \phi:C^*(E)/J\rightarrow\rho_{M,t}(C^*(E))/\rho_{M,t}(J)=
   \rho_{M,t}(C^*(E))/\K(\H_M) $$ 
given by $\phi(x+J)=\pi(\rho_{M,t}(x))$. This follows from 
the gauge-invariant uniqueness theorem \cite[Theorem 2.1]{bhrs}. 
Indeed, by \cite[Corollary 3.5]{bhrs} $C^*(E)/J$ is naturally 
isomorphic to $C^*((E/K_M)\setminus\beta(B_M))$, and the gauge action 
on $C^*(E)/J$ is inherited from the gauge action $\gamma$ on 
$C^*(E)$, since $J$ is gauge-invariant. We also need a matching 
action on $\pi(\rho_{M,t}(C^*(E)))$. For 
$z\in\T$ let $U_z$ be a unitary operator on $\H_M$ such that 
$U_z(\xi_\alpha)=z^{|\alpha|}\xi_\alpha$ for $\alpha\in 
A_M\cup L_M^0$. Since $U_z(\rho_{M,t}(p_v))U_z^*=
\rho_{M,t}(p_v)$ for $v\in E^0$ and 
$$ U_z(\rho_{M,t}(s_e))U_z^*=\left\{ \begin{array}{ll} 
   \rho_{M,t}(s_e) & \text{if } e\in L_M^1 \\ 
   z\rho_{M,t}(s_e) & \text{otherwise} \end{array} \right. $$ 
for $e\in E^1$, $\text{Ad}\,U_z$ is an automorphism of 
$\rho_{M,t}(C^*(E))$, and it induces an automorphism of  
$\pi(\rho_{M,t}(C^*(E)))$. Therefore we may define an  
action $\theta$ of $\T$ on $\pi(\rho_{M,t}(C^*(E)))$ by 
$$ \theta_z(\pi(\rho_{M,t}(x)))=\pi(U_z\rho_{M,t}(x)U_z^*). $$ 
For all $x\in C^*(E)$ and $z\in\T$ we have 
$\theta_z(\phi(x+J))=\phi(\gamma_z(x)+J)$, 
since this identity holds on the generators $\{p_v,s_e\}$ 
of $C^*(E)$. We must still show that $\phi$ does not kill 
any of the generating projections of $C^*(E)/J\cong 
C^*((E/K_M)\setminus\beta(B_M))$. We set $K'=E^0\setminus K_M=
M\setminus K_M$, $B'=\KMbad\setminus B_M$. Since $((E/K_M)\setminus
\beta(B_M))^0=K'\cup\{\beta(w):w\in B'\}$ it suffices to 
show that $\rho_{M,t}(p_{v})\not\in\K(\H_M)$ for $v\in K'$ 
and $\rho_{M,t}(p_w-p_{w,K_M})\not\in\K(\H_M)$ 
for $w\in B'$. That is, we must prove that 
$\rho_{M,t}(p_v)$ and $\rho_{M,t}(p_w-p_{w,K_M})$ are 
infinite dimensional for $v\in K'$ and $w\in B'$. For $v\in K'$ this 
fact is simply contained in the definition of $K_M$. Since $w\in B'$ 
there are infinitely many edges $e\in E^1$ such that 
$s(e)=w$ and $r(e)\in K_M\setminus\Omega(M)$, and consequently 
$\rho_{M,t}(p_w-p_{w,K_M})$ is infinite dimensional. Thus the 
hypothesis of the gauge-invariant uniqueness theorem is 
satisfied, and we may conclude that $\ker\rho_{M,t}\subseteq J$. 

Now let $J_2$ be a gauge-invariant ideal of $C^*(E)$ 
containing $\ker\rho_{M,t}$. We must show that $J\subseteq J_2$. 
It follows from part 1 that $J_{\Omega(M),\OMbad}\subseteq J_2$. 
Since $p_{r(L_M)}-(1/t)s_{L_M}\in\ker\rho_{M,t}\subseteq J_2$ 
and $J_2$ is gauge-invariant, also $p_{r(L_M)}=\Gamma(
p_{r(L_M)}-(1/t)s_{L_M})\in J_2$ and hence  
$\{p_v:v\in L_M^0\}\subseteq J_2$. Now if 
$v\in K_M\setminus(\Omega(M)\cup L_M^0)$ then $p_v\in J_2$ by 
Lemma \ref{projectionsinKM}. Consequently $I_{K_M}\subseteq J_2$. 
If $v\in B_M$ then the finite sum $\sum_{s(e)=v,r(e)\in K_M
\setminus\Omega(M)}s_e s_e^*$ belongs to $I_{K_M}$, and 
$p_v-p_{v,\Omega(M)}$ belongs to $J_{\Omega(M),\OMbad}$. Thus 
$$ p_v-p_{v,K_M}=\left( p_v-p_{v,\Omega(M)}\right)+
   \sum_{s(e)=v,r(e)\in K_M\setminus\Omega(M)}s_e s_e^* $$ 
belongs to $J_{\Omega(M),\OMbad}+I_{K_M}
\subseteq J_2$, and consequently $J=J_{K_M,B_M}
\subseteq J_2$, as required. 
\end{proof} 

In particular, for each $M\in{\mathcal M}_\tau(E)$ and $t\in\T$ we have 
$$ J_{\Omega(M),\OMbad}\subset\ker\rho_{M,t}\subset J_{K_M,B_M}. $$ 

From Lemmas \ref{rhoirreducible} and \ref{sandwich} 
we deduce the following. 

\begin{cor}\label{compacts} 
Let $E$ be a directed graph, $M\in{\mathcal M}_\tau(E)$ and 
$t\in\T$. Then 
$$ J_{K_M,B_M}=(\rho_{M,t})^{-1}(\K(\H_M)). $$  
\end{cor} 

In the following Lemma \ref{loopquotient} we find explicit 
generators for the ideals $\ker\rho_{M,t}$. The lemma also 
implies that for each $M\in{\mathcal M}_\tau(E)$ the family 
$\{\ker\rho_{M,t}:t\in\T\}$ imbeds topologically as a circle 
into the primitive ideal space of $C^*(E)$. 

\begin{lemma}\label{loopquotient} 
Let $E$ be a directed graph, $M\in{\mathcal M}_\tau(E)$, $v=r(L_M)$,   
and $\pi:J_{K_M,B_M}\rightarrow J_{K_M,B_M}/J_{\Omega(M),\OMbad}$ 
be the canonical surjection. Then the following hold. 
\begin{enumerate} 
\item The hereditary $C^*$-subalgebra 
$\pi(p_v)(J_{K_M,B_M}/J_{\Omega(M),\OMbad})\pi(p_v)$ is a full corner 
in the quotient $J_{K_M,B_M}/J_{\Omega(M),\OMbad}$. It is 
generated by $\pi(s_{L_M})$ and isomorphic to $C(\T)$. Hence the 
quotient $J_{K_M,B_M}/J_{\Omega(M),\OMbad}$ is Morita equivalent to $C(\T)$. 
\item For $t\in\T$ the ideal $\ker\rho_{M,t}$ of $C^*(E)$ 
is generated  by $J_{\Omega(M),\OMbad}$ and 
$s_{L_M}-tp_v$. We have $\ker\rho_{M,t}=\ker\rho_{L_M',t}$ for any 
cyclic permutation $L_M'$ of $L_M$,  and $\ker\rho_{M,t}
\neq\ker\rho_{M,z}$ if $t\neq z\in\T$. 
\end{enumerate} 
\end{lemma} 
\begin{proof} 
Ad 1. 
Let $J=\overline{\text{span}}\{J_{K_M,B_M}p_vJ_{K_M,B_M}+
J_{\Omega(M),\OMbad}\}$, an ideal of $C^*(E)$ contained in 
$J_{K_M,B_M}$. By Lemma \ref{projectionsinKM}, 
$\{p_w:w\in K_M\setminus(\Omega(M)\cup L_M^0)\}\subseteq J$, 
and clearly $\{p_w:w\in\Omega(M)
\cup L_M^0\}\subseteq J$. Thus $I_{K_M}\subseteq J$. Also, if 
$w\in B_M$ then $p_w-p_{w,K_M}\in I_{K_M}+J_{\Omega(M),\OMbad}$ 
(cf. the argument at the end of the proof of Lemma \ref{sandwich}), 
and thus $p_w-p_{w,K_M}\in J$. Consequently $J_{K_M,B_M}=J$ 
and hence $\pi(p_v)(J_{K_M,B_M}/J_{\Omega(M),\OMbad})\pi(p_v)$ 
is a full corner in $J_{K_M,B_M}/J_{\Omega(M),\OMbad}$. 

If $\mu,\nu\in E^*$ then $p_v s_\mu s_\nu^* p_v
\not\in J_{\Omega(M),\OMbad}$ if and only if 
$s(\mu)=s(\nu)=v$ and $r(\mu)=r(\nu)\in L_M^0$. Thus 
$\pi(p_v)(J_{K_M,B_M}/J_{\Omega(M),\OMbad})\pi(p_v)$  is 
generated as a $C^*$-algebra by $\pi(s_{L_M})$. If 
$z\in\T$ then $J_{\Omega(M),\OMbad}\subseteq\ker\rho_{M,z}$ 
and, consequently, $\rho_{M,z}$ induces a representation $\tilde{\rho}_{M,z}$ 
of $J_{K_M,B_M}/J_{\Omega(M),\OMbad}$. Since $\tilde{\rho}_{M,z}
(\pi(s_{L_M}))$ equals $z$-multiple of a rank one projection, 
it follows that the spectrum 
of the partial unitary $\pi(s_{L_M})$ contains the entire unit circle. 
Thus, the corner $\pi(p_v)(J_{K_M,B_M}/J_{\Omega(M),\OMbad})
\pi(p_v)$  is isomorphic to $C(\T)$. 

Ad 2. 
Fix $t\in\T$ and let $J'$ be the ideal of $C^*(E)$ 
generated by $J_{\Omega(M),\OMbad}$ and $s_{L_M}-tp_v$. 
Since $\rho_{M,t}(s_{L_M}-tp_v)=0$, Lemma \ref{sandwich} implies that 
$$ J_{\Omega(M),\OMbad}\subseteq J'\subseteq\ker\rho_{M,t}
   \subseteq J_{K_M,B_M}. $$ 
We have shown, above, that $\pi(J_{K_M,B_M}p_v
J_{K_M,B_M})=\pi(J_{K_M,B_M})$ and that $\pi(p_vJ_{K_M,B_M}
p_v)=C^*(\pi(s_{L_M}))\cong C(\T)$. Hence  
$\pi(p_vJ'p_v)=\{g(\pi(s_{L_M})):g\in C(\T),g(t)=0\}=
\pi(p_v(\ker\rho_{M,t})p_v)$. Thus 
\begin{align*} 
\pi(J') & =  \pi((J_{K_M,B_M}p_vJ_{K_M,B_M})
   J'(J_{K_M,B_M}p_vJ_{K_M,B_M})) \\ 
 & =  \pi(J_{K_M,B_M}p_v(\ker\rho_{M,t})p_vJ_{K_M,B_M}) \\  
 & =  \pi((J_{K_M,B_M}p_vJ_{K_M,B_M})
   (\ker\rho_{M,t})(J_{K_M,B_M}p_vJ_{K_M,B_M})) \\ 
 & =  \pi(\ker\rho_{M,t}).  
\end{align*} 
It follows that $J'=\ker\rho_{M,t}$ and consequently the ideal 
$\ker\rho_{M,t}$ is generated by $J_{\Omega(M),\OMbad}$ and 
$s_{L_M}-tp_v$. This immediately implies that $\ker\rho_{M,t}=
\ker\rho_{L_M',t}$ for any cyclic permutation $L_M'$ of $L_M$. 
Finally, if $t\neq z$ then $s_{L_M}-tp_v\in\ker\rho_{M,t}
\setminus\ker\rho_{M,z}$ and hence $\ker\rho_{M,t}\neq\ker\rho_{M,z}$. 
\end{proof} 

For $M\in{\mathcal M}_\tau(E)$ and $t\in\T$ we denote by 
$R_{M,t}$ the closed two-sided ideal of $C^*(E)$ generated by 
$J_{\Omega(M),\OMbad}$ and $s_{L_M}-tp_{r(L_M)}$. By Lemma 
\ref{loopquotient}  
$$ R_{M,t}=\ker\rho_{M,t}. $$ 
Since $\gamma_{\overline{t}}(R_{M,1})=R_{M,t}$, 
these ideals are not gauge-invariant. 
We will show in Theorem \ref{otherprimitive}  
that each non gauge-invariant primitive ideal of $C^*(E)$ 
is of the form $R_{M,t}$. To this end we still need the  
following simple lemma. 

\begin{lemma}\label{condJprimitive} 
Let $E$ be a directed graph. If $J\neq0$ is 
a primitive ideal of $C^*(E)$ such that 
$p_v\not\in J$ for all $v\in E^0$, then $E^0\in{\mathcal M}_\tau(E)$ 
and there is a $t\in\T$ such that $J=R_{E^0,t}$. 
\end{lemma}
\begin{proof} 
Since $p_v\not\in J$ for all $v\in E^0$, 
it follows from \cite[Lemma 4.1]{bhrs} that $E^0$ is a maximal tail.  
If all loops in $E$ had exits then there existed a $v\in E^0$ 
such that $p_v\in J$, by the Cuntz-Krieger uniqueness theorem  
\cite[Theorem 2]{flr}. Thus $E^0\in{\mathcal M}_\tau(E)$. 
To simplify the notation, in the remaining part of this proof 
we write $M=E^0$. 

Let $\rho$ be an irreducible representation of $C^*(E)$ with 
kernel $J$. Since $J$ does not contain any projections $p_v$, 
$v\in E^0$, $J$ does not contain the ideal $I_{L_M^0}$. Thus 
the restriction of $\rho$ to $I_{L_M^0}$ must be irreducible. 
By Lemma \ref{loopquotient} there exists a $t\in\T$ such 
that the restrictions of $\rho$ and $\rho_{M,t}$ to $I_{L_M^0}$ 
coincide. Hence $\rho=\rho_{M,t}$ and consequently 
$J=\ker\rho=\ker\rho_{M,t}=R_{M,t}$. 
\end{proof} 

We set $\text{Prim}_\tau(C^*(E))=\text{Prim}(C^*(E))
\setminus\text{Prim}_\gamma(C^*(E))$, the collection of primitive 
ideals of $C^*(E)$ which are not invariant under 
the gauge action $\gamma$. 

\begin{theo}\label{otherprimitive} 
Let $E$ be a directed graph. The map 
$$ {\mathcal M}_\tau(E)\times\T\longrightarrow\text{Prim}_\tau(C^*(E)) $$ 
given by 
$$ (M,t)\longmapsto R_{M,t} $$  
is a bijection. 
\end{theo} 
\begin{proof} 
The map is well-defined by Lemmas \ref{rhoirreducible} 
and \ref{loopquotient}. 

Firstly, we show that the map is injective. That is, we must show that 
the ideals $\{R_{M,t}:M\in{\mathcal M}_\tau(E),t\in\T\}$ are distinct. 
Indeed, if $R_{M,t}=R_{N,z}$ then $J_{\Omega(M),\OMbad}=J_{\Omega(N),
\ONbad}$ by Lemma \ref{sandwich}. Thus $\Omega(M)=\Omega(N)$ by 
\cite[Lemma 3.7]{bhrs}, and consequently $M=E^0\setminus\Omega(M)=
E^0\setminus\Omega(N)=N$. It then follows from Lemma \ref{loopquotient}
that $t=z$. 

Secondly, we show that the map is surjective. 
Let $J\in\text{Prim}_\tau(C^*(E))$. We set $K=\{v\in E^0:p_v\in J\}$ and 
$B=\{x\in \Kbad:p_x-p_{x,K}\in J\}$. Then $J_{K,B}$ is a proper ideal of $J$ 
and hence $J/J_{K,B}$ is a non-zero primitive ideal of $C^*(E)/J_{K,B}$. 
By \cite[Corollary 3.5]{bhrs} we have $C^*(E)/J_{K,B}\cong C^*(F)$, 
with $F=(E/K)\setminus\beta(B)$. We denote the canonical generating 
Cuntz-Krieger $F$-family by $\{q_w,u_f\}$ (cf. \cite[Proposition 3.4]{bhrs}). 
By \cite[Lemma 3.7]{bhrs} the ideal $J/J_{K,B}$ does not contain any 
projections $q_w$, $w\in F^0$. Now applying Lemma \ref{condJprimitive} 
to the ideal $J/J_{K,B}$ of $C^*(F)$ we see that $F^0\in
{\mathcal M}_\tau(F)$ and, using also Lemma \ref{loopquotient}, that there is 
$t\in\T$ such that $J/J_{K,B}$ is generated as an ideal of $C^*(F)$ by 
$u_{L_{F^0}}-tq_{r(L_{F^0})}$. 

Let $M=E^0\setminus K$. $M$ is a maximal tail in $E$ by 
\cite[Lemma 4.1]{bhrs}. If $L_{F^0}=(f_1,\ldots,f_k)$ then all $f_i$ 
must come from edges in $E^1$. Clearly the loop $L_{F^0}$ 
has no exits in $M$. Thus $M\in{\mathcal M}_\tau(E)$ and 
$L_{F^0}$ is a cyclic permutation of $L_M$ by Lemma 
\ref{looptail}. Since $M$ is a maximal tail in $E$ we have $K=E^0\setminus M
=\Omega(M)$. We also have $B=\OMbad$. Indeed, otherwise the graph $F$ 
would contain a sink  $\beta(v)$, $v\in\OMbad\setminus B$, contradicting 
the fact that $F^0$ belongs to ${\mathcal M}_\tau(F)$. Consequently 
the ideal $J$ contains $J_{\Omega(M),\OMbad}=J_{K,B}$. Since 
$u_{L_{F^0}}-tq_{r(L_{F^0})}$ belongs to the quotient 
$J/J_{\Omega(M),\OMbad}$, it now follows that $s_{L_M}-tp_{r(L_M)}$ 
belongs to $J$. By Lemma \ref{loopquotient} we have $R_{M,t}\subseteq J$. 
As both $J/J_{\Omega(M),\OMbad}$ and $R_{M,t}/J_{\Omega(M),\OMbad}$ 
are generated by $u_{L_{F^0}}-tq_{r(L_{F^0})}$, it follows that $J=R_{M,t}$. 
\end{proof} 

Combining Theorem \ref{otherprimitive} with \cite[Theorem 4.7]{bhrs} 
we obtain a complete list of primitive ideals of $C^*(E)$ for an 
arbitrary countable graph $E$.  

\begin{cor}\label{primitiveideals} 
For a countable directed graph $E$ the map 
$$ {\mathcal M}_\gamma(E)\cup BV(E)\cup({\mathcal M}_\tau(E)
   \times\T)\longrightarrow\text{Prim}(C^*(E)) $$ 
given by 
\begin{align*} 
{\mathcal M}_{\gamma}(E)\ni M & \longmapsto  J_{\Omega(M),\OMbad} \\ 
BV(E)\ni v & \longmapsto  J_{\Omega(v),\Ovbad\setminus\{v\}} \\ 
{\mathcal M}_\tau(E)\times\T\ni(N,t) & \longmapsto  R_{N,t} 
\end{align*} 
is a bijection. 
\end{cor}  

If $E$ is row-finite then $BV(E)=\emptyset$ and 
$\Kbad=\emptyset$ for every saturated hereditary $K\subseteq E^0$. 
Consequently, for such graphs we have the following 
simpler description of the primitive ideals. 

\begin{cor}\label{primitiverowfinite} 
For a countable, row-finite directed graph $E$ the map 
$$ {\mathcal M}_\gamma(E)\cup({\mathcal M}_\tau(E)\times\T)
   \longrightarrow\text{Prim}(C^*(E)) $$ 
given by 
\begin{align*} 
{\mathcal M}_{\gamma}(E)\ni M & \longmapsto  I_{\Omega(M)} \\ 
{\mathcal M}_\tau(E)\times\T\ni(N,t) & \longmapsto  R_{N,t} 
\end{align*} 
is a bijection. 
\end{cor} 

\section{The hull-kernel topology} 

$\text{Prim}(C^*(E))$ is a topological space with the 
hull-kernel topology determined by the closure operation 
$$ \overline{Z}=\{J\in\text{Prim}(C^*(E)):\cap Z\subseteq J\}. $$ 
Our goal in this section is to describe this closure operation. 
Using the bijection of Corollary \ref{primitiveideals} we 
transport the hull-kernel topology from $\Prim(C^*(E))$ 
onto ${\mathcal M}_\gamma(E)\cup 
BV(E)\cup({\mathcal M}_\tau(E)\times\T)$. We begin with the 
following two simple lemmas. 

\begin{lemma}\label{rhoinclusion} 
Let $E$ be a directed graph and $M\neq N\in{\mathcal M}_\tau(E)$.  
If there exists a path from $L_{N}^0$ to $L_M^0$ then 
for all $t,z\in\T$ we have 
$$ R_{M,t}\subset J_{K_M,B_M}\subseteq 
   J_{\Omega(N),\ONbad}\subset R_{N,z}. $$  
\end{lemma} 
\begin{proof} 
By virtue of Lemma \ref{sandwich} it suffices to prove that 
$R_{M,t}\subseteq J_{\Omega(N),\ONbad}$. By Lemma \ref{loopquotient} 
this amounts to showing that $J_{\Omega(M),\OMbad}\subseteq 
J_{\Omega(N),\ONbad}$ and $s_{L_M}-tp_{r(L_M)}\in J_{\Omega(N),\ONbad}$. 
Indeed, since there is a path from $L_{N}^0$ to $L_M^0$ we have 
$\Omega(M)\subseteq\Omega(N)$ and $\OMbad\subseteq\Omega(N)\cup\ONbad$. 
Thus $J_{\Omega(M),\OMbad}\subseteq J_{\Omega(N),\ONbad}$ by 
\cite[Corollary 3.10]{bhrs}. Furthermore, Lemma \ref{looptail} implies 
that there is no path from $L_M^0$ to $L_{N}^0$. Hence both $p_{r(L_M)}$ 
and $s_{L_M}=s_{L_M}p_{r(L_M)}$ are in $J_{\Omega(N),\ONbad}$.  
\end{proof} 

\begin{lemma}\label{mtintersection} 
Let $E$ be a directed graph. If $Y\subseteq{\mathcal M}_\tau(E)$ 
then we have 
$$ \bigcap_{U\in Y}J_{\Omega(U),\OUbad}=J_{K,\Kbad} \;\;
   \text{ with } \;\; K=\bigcap_{U\in Y}\Omega(U). $$ 
\end{lemma} 
\begin{proof} 
By \cite[Proposition 3.9]{bhrs} we have $J=J_{K,B}$ with 
$K=\bigcap_{U\in Y}\Omega(U)$ and $B=(\bigcap_{U\in Y}\Omega(U)
\cup\OUbad)\cap\Kbad$. Fix $U_0\in Y$ and let $w\in\Kbad\setminus
\Omega(U_0)$. Then $w$ emits infinitely many edges into 
$K=\bigcap_{U\in Y}\Omega(U)\subseteq\Omega(U_0)$ and only finitely 
many edges outside $K$, hence also only finitely 
many edges outside $\Omega(U_0)$. Since there is a path 
from $w$ to the loop $L_{U_0}$, $w$ emits at least one edge 
into $E^0\setminus\Omega(U_0)$. Consequently $\Kbad\subseteq
\bigcap_{U\in Y}\Omega(U)\cup\OUbad$ and hence $B=\Kbad$.  
\end{proof} 

If $K$ is a hereditary saturated subset of $E^0$ then  
the set $\Kbad$ of vertices which emit infinitely many edges 
into $K$ and finitely many edges into its complement 
affects the ideal structure of $C^*(E)$ and hence it affects 
$\Prim(C^*(E))$. To describe the topology of $\Prim(C^*(E))$ it 
is also important to consider the set of those vertices which 
emit infinitely many edges into $K$ and none into its complement. 
We call this set $K_\infty^\emptyset$. More formally, we define 
$$ K_\infty^\emptyset:=\{v\in E^0\setminus K:|s^{-1}(v)|=\infty \text{ and } 
   r(e)\in K \text{ for all } e \text{ with } s(e)=v\}. $$ 
If $M$ is a maximal tail then $\M0$ is either empty or consists 
of exactly one element by (MT3). In the latter case, if 
$\M0=\{w\}$ then $M$ consists of all those vertices $u\in E^0$ 
for which there exists a path from $u$ to $w$. 

\begin{lemma}\label{3inclusions} 
Let $E$ be a directed graph. Let $M\in{\mathcal M}_\gamma(E)$, 
$v\in BV(E)$, $N\in{\mathcal M}_\tau(E)$, and $t\in\T$. 
If $Y\subseteq{\mathcal M}_\tau(E)$ and $K=\bigcap_{U\in Y}\Omega(U)$ 
then the following hold. 
\begin{enumerate} 
\item $J_{K,\Kbad}\subseteq J_{\Omega(M),\OMbad}$ if and 
only if $M\subseteq\bigcup Y$ and 
$s^{-1}(\M0)\cap r^{-1}(K)$ is finite. 
\item $J_{K,\Kbad}\subseteq J_{\Omega(v),\Ovbad\setminus\{v\}}$ 
if and only if $v\in \bigcup Y$ and 
$s^{-1}(v)\cap r^{-1}(K)$ is finite. 
\item $J_{K,\Kbad}\subseteq R_{N,t}$ 
if and only if $N\subseteq\bigcup Y$.  
\end{enumerate} 
\end{lemma} 
\begin{proof} 
Ad 1. 
By \cite[Corollary 3.10]{bhrs}, $J_{K,\Kbad}
\subseteq J_{\Omega(M),\OMbad}$ if and 
only if $K\subseteq\Omega(M)$ and $\Kbad\subseteq\Omega(M)
\cup\OMbad$. Clearly, $K\subseteq\Omega(M)$ 
if and only if $M\subseteq\bigcup Y$. Now assuming  
$M\subseteq\bigcup Y$ we automatically have $\Kbad\setminus 
\M0\subseteq \Omega(M)\cup\OMbad$. On the other hand, if $w\in 
\M0$ then $w\not\in\Omega(M)\cup\OMbad$. Thus we must have 
$w\not\in\Kbad$ and this can only happen if 
$s^{-1}(\M0)\cap r^{-1}(K)$ is finite. 

Ad 2. 
By \cite[Corollary 3.10]{bhrs}, 
$J_{K,\Kbad}\subseteq J_{\Omega(v),\Ovbad
\setminus\{v\}}$ if and only if $K\subseteq\Omega(v)$ and 
$\Kbad\subseteq \Omega(v)\cup(\Ovbad\setminus\{v\})$. 
$K\subseteq\Omega(v)$ if and only if $v\in\bigcup Y$. 
Assuming $v\in\bigcup Y$ we have $(\Kbad\setminus\{v\})
\subseteq \Omega(v)\cup(\Ovbad
\setminus\{v\})$. Since $v\not\in\Omega(v)\cup(\Ovbad
\setminus\{v\})$ we must have $v\not\in \Kbad$, 
which can only happen if $s^{-1}(v)\cap r^{-1}(K)$ is finite. 

Ad 3. 
Since $J_{K,\Kbad}$ is gauge-invariant it follows from 
Lemma \ref{sandwich} that $J_{K,\Kbad}\subseteq 
R_{N,t}$ if and only if $J_{K,\Kbad}\subseteq J_{\Omega(N),\ONbad}$. 
For $K\subseteq\Omega(N)$ we must have $N\subseteq\bigcup Y$, and under 
this assumption we automatically have $\Kbad\subseteq \Omega(N)
\cup\ONbad$, since $\Omega(N)_\infty^\emptyset=\emptyset$. 
\end{proof} 

If $Y\subseteq{\mathcal M}_\tau(E)$ then it is convenient to 
consider two special subsets of $Y$, $Y_{\text{\rm min}}$ and $Y_\infty$, 
defined as follows. 
\begin{align*} 
Y_{\text{\rm min}} & :=  \{U\in Y: \text{ for all } U'\in Y,\, U'\neq U\, 
\text{ there is no path from } L_U^0 \text{ to } L_{U'}^0\}, \\ 
Y_\infty & :=  \{U\in Y: \text{ for all } V\in Y_{\text{\rm min}} 
\text{ there is no path from } L_U^0 \text{ to } L_V^0\}. 
\end{align*} 
We call $Y_{\text{\rm min}}$ the set of minimal elements of $Y$. 
We are now ready to describe the closure operation in $\Prim(C^*(E))$. 

\begin{theo}\label{9closures} 
Let $E$ be a countable directed graph. Let $X\subseteq
{\mathcal M}_\gamma(E)$,  $W\subseteq BV(E)$, $Y\subseteq
{\mathcal M}_\tau(E)$, and let $D(U)\subseteq\T$ for each $U\in Y$. 
If $M\in{\mathcal M}_\gamma(E)$, $v\in BV(E)$, 
$N\in{\mathcal M}_\tau(E)$, and $z\in\T$, 
then the following hold. 
\begin{enumerate} 
\item $M\in\overline{X}$ if and only if either 
\begin{description} 
\item{(i)} $M\in X$ or 
\item{(ii)} $M\subseteq\bigcup X$ and $s^{-1}(\M0)\cap 
r^{-1}(\bigcap_{U\in X}\Omega(U))$ is finite. 
\end{description} 
\item $v\in\overline{X}$ if and only if $v\in
\bigcup X$ and $s^{-1}(v)\cap r^{-1}(\bigcap_{U\in X}\Omega(U))$ 
is finite. 
\item $(N,z)\in\overline{X}$ if and only if 
$N\subseteq\bigcup X$.   
\item $M\in\overline{W}$ if and only if $M\subseteq
E^0\setminus\bigcap_{w\in W}\Omega(w)$ and 
$s^{-1}(\M0)\cap r^{-1}(\bigcap_{w\in W}\Omega(w))$ 
is finite. 
\item $v\in\overline{W}$ if and only if either 
\begin{description} 
\item{(i)} $v\in W$ or 
\item{(ii)} $v\in E^0\setminus\bigcap_{w\in W}\Omega(w)$ and 
$s^{-1}(v)\cap r^{-1}(\bigcap_{w\in W}\Omega(w))$ is finite. 
\end{description} 
\item $(N,z)\in\overline{W}$ if and only if 
$N\subseteq E^0\setminus\bigcap_{w\in W}\Omega(w)$.   
\item $M$ is in the closure of $\{(U,t):U\in Y,\:t\in D(U)\}$ 
if and only if either 
\begin{description} 
\item{(i)} $M\subseteq\bigcup Y_\infty$ and 
$s^{-1}(\M0)\cap r^{-1}(\bigcap_{U\in Y_\infty}\Omega(U))$ 
is finite or  
\item{(ii)} $M\subseteq\bigcup Y_{\text{\rm min}}$ and 
$s^{-1}(\M0)\cap r^{-1}(\bigcap_{U\in Y_{\text{\rm min}}}\Omega(U))$ 
is finite. 
\end{description} 
\item $v$ is in the closure of 
$\{(U,t):U\in Y,\:t\in D(U)\}$ if and only if either 
\begin{description} 
\item{(i)} $v\in\bigcup Y_\infty$ and $s^{-1}(v)\cap r^{-1}
(\bigcap_{U\in Y_\infty}\Omega(U))$ is finite or 
\item{(ii)} $v\in\bigcup Y_{\text{\rm min}}$ and $s^{-1}(v)
\cap r^{-1}(\bigcap_{U\in Y_{\text{\rm min}}}\Omega(U))$ is finite. 
\end{description} 
\item $(N,z)$ is in the closure 
of $\{(U,t):U\in Y,\:t\in D(U)\}$ if and only if one of the following 
three conditions holds. 
\begin{description} 
\item{(i)} $N\subseteq\bigcup Y_\infty$.  
\item{(ii)} $N\not\in Y_{\text{\rm min}}$ and $N\subseteq
\bigcup Y_{\text{\rm min}}$.  
\item{(iii)} $N\in Y_{\text{\rm min}}$ and $z\in\overline{D(N)}$. 
\end{description} 
\end{enumerate} 
\end{theo} 
\begin{proof} 
Throughout the proof of cases 1--6 we denote 
\begin{align*} 
K=\bigcap_{U\in X}\Omega(U),\;\;\;\;\; & 
B=\left(\bigcap_{U\in X}\Omega(U)\cup
            \OUbad\right)\cap\Kbad, \\ 
K'=\bigcap_{w\in W}\Omega(w),\;\;\;\;\; & 
B'=\left(\bigcap_{w\in W}\Omega(w)\cup(\Owbad
         \setminus\{w\})\right)\cap\Kpbad. 
\end{align*} 

Ad 1. 
It suffices to consider the case $M\not\in X$. 
We have $\bigcap_{U\in X}J_{\Omega(U),\OUbad}=J_{K,B}$ 
by \cite[Proposition 3.9]{bhrs}. Thus $M\in\overline{X}$ if 
and only if $J_{K,B}\subseteq J_{\Omega(M),\OMbad}$. 
By \cite[Corollary 3.10]{bhrs} this is equivalent to 
$K\subseteq\Omega(M)$ and $B\subseteq \Omega(M)
\cup\OMbad$. Clearly, $K\subseteq\Omega(M)$ 
if and only if $M\subseteq\bigcup X$. Assuming  
$M\subseteq\bigcup X$ we automatically have $(B\setminus 
\M0)\subseteq \Omega(M)\cup\OMbad$. On the other hand, if 
$w\in \M0$ then $w\not\in \Omega(M)\cup\OMbad$. Thus  
we must have $w\not\in B$, which can only happen in one of the 
following three cases: (i) $\{r(e):e\in E^1,\:s(e)=w\}
\subseteq\bigcap_{U\in X}\Omega(U)$, (ii) there is a $U\in X$ 
such that $\{r(e):e\in E^1,\:s(e)=w\}\subseteq\Omega(U)$, 
(iii) $s^{-1}(w)\cap r^{-1}(\bigcap_{U\in X}\Omega(U))$ is finite. 
Case (i) reduces to case (ii) since we assumed that 
$M\subseteq\bigcup X$. In case (ii) we have $w\in\Omega(U)_\infty^0$ 
and hence $M=U\in X$, contrary to the assumption. Therefore 
only case (iii) remains, and the claim is proved. 

Ad 2. 
Similarly as in case 1 above, $v\in\overline{X}$ if and only 
if $J_{K,B}\subseteq J_{\Omega(v),\Ovbad\setminus\{v\}}$, and this   
is equivalent to $K\subseteq\Omega(v)$ and $B\subseteq\Omega(v)
\cup(\Ovbad\setminus\{v\})$. Clearly, $K\subseteq\Omega(v)$ 
if and only if $v\in\bigcup X$. Assuming  
$v\in\bigcup X$ we automatically have $(B\setminus 
\{v\})\subseteq \Omega(v)\cup(\Ovbad
\setminus\{v\})$. On the other hand, $v\not\in
\Omega(v)\cup(\Ovbad\setminus\{v\})$. 
Thus we must have $v\not\in B$, which can only happen if 
$s^{-1}(v)\cap r^{-1}(\bigcap_{U\in X}\Omega(U))$ is finite. 

Ad 3. 
Since $\bigcap_{U\in X}J_{\Omega(U),\OUbad}$ 
is gauge-invariant, by Lemma \ref{sandwich}  
this ideal is contained in $R_{N,z}$ 
if and only if it is already contained in 
$J_{\Omega(N),\ONbad}$. Since $N\not\in X$ and 
$\Omega(N)_\infty^\emptyset=\emptyset$ the claim is proved 
similarly to the first part of case 1 above. 

Ad  4. 
We have $\bigcap_{w\in W}
J_{\Omega(w),\Owbad\setminus\{w\}}=
J_{K',B'}$ . Thus, by \cite[Corollary 3.10]{bhrs}, 
$M\in\overline{W}$ if and only if $K'\subseteq\Omega(M)$ 
and $B'\subseteq\Omega(M)\cup \OMbad$. 
$K'\subseteq\Omega(M)$ is equivalent to $M\subseteq
E^0\setminus\bigcap_{w\in W}\Omega(w)$. Assuming this,   
it suffices to find a condition for $v\not\in B'$ 
(similarly to the argument from case 1 above). It is easy to 
see that $v\not\in B'$ occurs precisely when $s^{-1}
(\M0)\cap r^{-1}(\bigcap_{w\in W}\Omega(w))$ is finite. 

Ad 5. Obviously, if $v\in W$ then $v\in\overline{W}$. Thus 
we may assume that $v\not\in W$. Similarly to 
the above, $v\in\overline{W}$ if and only if $K'\subseteq
\Omega(v)$ and $B'\subseteq\Omega(v)\cup(\Ovbad
\setminus\{v\})$. $K'\subseteq\Omega(v)$ is equivalent to 
$v\in E^0\setminus\bigcap_{w\in W}\Omega(w)$. Since 
$B'\setminus\{v\}\subseteq\Omega(v)\cup(\Ovbad\setminus\{v\})$ 
and $v\not\in\Omega(v)\cup(\Ovbad\setminus\{v\})$, it suffices 
to find a condition for $v\not\in B'$. But this is equivalent 
to $s^{-1}(v)\cap r^{-1}(\bigcap_{w\in W}\Omega(w))$ being finite. 

Ad 6. 
The proof is similar to the case 3 above. Indeed, 
$(N,z)\in\overline{W}$ if and only if $J_{K',B'}
\subseteq R_{N,z}$. Since $J_{K',B'}$ 
is gauge-invariant this is equivalent to $J_{K',B'}
\subseteq J_{\Omega(N),\ONbad}$ by Lemma 
\ref{sandwich}, and by \cite[Corollary 3.10]{bhrs} this 
happens if and only if $K'\subseteq\Omega(N)$ and $B'
\subseteq\Omega(N)\cup\ONbad$. The latter is 
automatically satisfied and the former is equivalent 
to $N\subseteq E^0\setminus\bigcap_{w\in W}\Omega(w)$.  

Ad 7--9.  
The following observations are used in the proofs of 
cases 7, 8 and 9. 

The closure of $\{(U,t):U\in Y,\:t\in D(U)\}$ coincides 
with the union of the closure of $\{(U,t):U\in Y_\infty,
\:t\in D(U)\}$ and the closure of $\{(U,t):U\in Y
\setminus Y_\infty, \:t\in D(U)\}$. Thus it suffices 
to find these two closures.  

In order to determine the closure of $\{(U,t):U\in Y_\infty,
\:t\in D(U)\}$ we observe that 
\begin{equation}\label{Ymc} 
\bigcap_{U\in Y_\infty}\bigcap_{t\in D(U)}R_{U,t}=
\bigcap_{U\in Y_\infty}J_{\Omega(U),\OUbad}.   
\end{equation}  
Indeed, if $U'\in Y_\infty$ then there exists a $U\in Y_\infty$ 
different from $U'$ such that there is a path from 
$L_{U'}^0$ to $L_U^0$. By Lemma \ref{rhoinclusion} we have 
$R_{U',t_1}\cap R_{U,t_2}=J_{\Omega(U'),\OUpbad}\cap R_{U,t_2}$ 
and hence in the LHS of (\ref{Ymc}) 
we may replace each $R_{U,t}$ by $J_{\Omega(U),
\OUbad}$. Therefore, a primitive ideal 
belongs to the closure of $\{(U,t):U\in Y_\infty,
\:t\in D(U)\}$ if and only if the conditions of Lemma 
\ref{3inclusions} are satisfied (with $Y_\infty$ instead of $Y$).  

Similarly, in order to determine the closure of $\{(U,t):U\in 
Y\setminus Y_\infty,\:t\in D(U)\}$ we observe that 
\begin{equation}\label{Ym}  
\bigcap_{U\in Y\setminus Y_\infty}\bigcap_{t\in D(U)}
R_{U,t}=\bigcap_{U\in Y_{\text{\rm min}}}\bigcap_{t\in D(U)}R_{U,t}  
\end{equation} 
by virtue of Lemma \ref{rhoinclusion}.  
By Lemmas  \ref{mtintersection} and \ref{sandwich} we get 
$$ J_{K'', \Kppbad}=\bigcap_{U\in Y_{\text{\rm min}}}J_{\Omega(U),\OUbad}
\subseteq\bigcap_{U\in Y_{\text{\rm min}}} \bigcap_{t\in D(U)}R_{U,t}, $$ 
where $K''=\bigcap_{U\in Y_{\text{\rm min}}}\Omega(U)$.   
Then a primitive ideal $J$ belongs to the closure of $\{(U,t):U\in 
Y\setminus Y_\infty,\:t\in D(U)\}$ if and only if $J$ contains the 
RHS of (\ref{Ym}), and for this it is necessary that $J\supseteq  
J_{K'',\Kppbad}$. Thus it is useful to look 
at the quotient $C^*(E)/J_{K'',\Kppbad}$. 
By \cite[Corollary 3.5]{bhrs} we have $C^*(E)/J_{K'',\Kppbad}
\cong C^*(F)$, where $F=(E/K'')\setminus\beta(\Kppbad)$ is the 
subgraph of $E$ such that $F^0=E^0\setminus K''$ and $F^1=\{e\in E^1:
r(e)\not\in K''\}$.  

Ad 7. 
$J_{\Omega(M),\OMbad}$ contains  
$J_{K'',\Kppbad}$ if and only if $M\subseteq\bigcup Y_{\text{\rm min}}$ 
and $s^{-1}(\M0)\cap r^{-1}(K'')$ is finite, by 
Lemma \ref{3inclusions}. Assume this holds. Then 
$J_{\Omega(M),\OMbad}/J_{K'',\Kppbad}$ 
is a gauge-invariant primitive ideal of $C^*(F)$ and hence 
contains all projections corresponding to $\{v\in L_U^0:
U\in Y_{\text{\rm min}}\}$, since the loops $L_U$, $U\in Y_{\text{\rm min}}$ 
have no exits in $F$. By Lemmas \ref{projectionsinKM} 
and \ref{sandwich} this implies that 
$R_{U,t}\subseteq J_{K_U,B_U}\subseteq  
J_{\Omega(M),\OMbad}$ for each $U\in Y_{\text{\rm min}}$, 
$t\in D(U)$, and thus $M$ is in the closure of 
$\{(U,t):U\in Y\setminus Y_\infty,\:t\in D(U)\}$. Consequently, 
$M$ is in the closure of $\{(U,t):U\in Y,\:t\in D(U)\}$ if 
and only if either (i) $M\subseteq\bigcup Y_\infty$ and 
$s^{-1}(\M0)\cap r^{-1}(\bigcap_{U\in Y_\infty}\Omega(U))$ 
is finite, or (ii) $M\subseteq\bigcup Y_{\text{\rm min}}$ and 
$s^{-1}(\M0)\cap r^{-1}(\bigcap_{U\in Y_{\text{\rm min}}}\Omega(U))$ is 
finite. 

Ad 8. 
This is proved by an argument very similar to case 7 above. 

Ad 9. 
$R_{N,z}$ contains $J_{K'',\Kppbad}$ if and only if $N\subseteq\bigcup 
Y_{\text{\rm min}}$, by Lemma \ref{3inclusions}. Assume this holds. 
If $N\not\in Y_{\text{\rm min}}$ then $R_{N,z}$ contains the RHS 
of (\ref{Ym}) by Lemma \ref{rhoinclusion}, since there exists 
a path from $L_{N}^0$ to at least one $L_U^0$, $U\in Y_{\text{\rm min}}$. 
Suppose $N\in Y_{\text{\rm min}}$. If $z\not\in\overline{D(N)}$ then 
let $g:\T\rightarrow\C$ be a continuous function such that 
$g|_{\overline{D(N)}}=0$ and $g(z)\neq0$. Then 
$g(s_{L_{N}})$ is in the RHS of (\ref{Ym}) but 
not in $R_{N,z}$. Thus we must have 
$z\in\overline{D(N)}$. In this case it follows from 
Lemma \ref{loopquotient} that $R_{N,z}$ 
contains the RHS of (\ref{Ym}). Consequently, $(N,z)$ 
is in the closure of $\{(U,t):U\in Y,\:t\in D(U)\}$ 
if and only if one of the following three conditions holds; 
(i) $N\subseteq\bigcup Y_\infty$, 
(ii) $N\not\in Y_{\text{\rm min}}$ and 
$N\subseteq\bigcup Y_{\text{\rm min}}$, 
(iii) $N\in Y_{\text{\rm min}}$ and $z\in\overline{D(N)}$. 
\end{proof} 

\begin{cor}\label{4closures} 
Let $E$ be a row-finite directed graph. Let $X\subseteq
{\mathcal M}_\gamma(E)$, $Y\subseteq{\mathcal M}_\tau(E)$, and let 
$D(U)\subseteq\T$ for each $U\in Y$. If $M\in{\mathcal M}_\gamma(E)$, 
$N\in{\mathcal M}_\tau(E)$, and $z\in\T$, then the following hold. 
\begin{enumerate} 
\item $M\in\overline{X}$ if and only if $M\subseteq\bigcup X$. 
\item $(N,z)\in\overline{X}$ if and only if $N\subseteq\bigcup X$.   
\item $M$ is in the closure of $\{(U,t):U\in Y,\:t\in D(U)\}$ 
if and only if $M\subseteq\bigcup Y$. 
\item $(N,z)$ is in the closure of $\{(U,t):U\in Y,\:t\in D(U)\}$ 
if and only if one of the following 
three conditions holds. 
\begin{description} 
\item{(i)} $N\subseteq\bigcup Y_\infty$.  
\item{(ii)} $N\not\in Y_{\text{\rm min}}$ and $N\subseteq
\bigcup Y_{\text{\rm min}}$.  
\item{(iii)} $N\in Y_{\text{\rm min}}$ and $z\in\overline{D(N)}$. 
\end{description} 
\end{enumerate} 
\end{cor} 

\section{Examples}

We illustrate the main results of this paper with the following three
examples. The discussion of gauge-invariant ideals of the  
algebras corresponding to the last two of them was carried out in 
\cite[Section 5]{bhrs}. Now we are in a position to give a complete 
description of their primitive ideal spaces. 

\begin{example}\label{example1} 
Let $E_1$ be the following graph, in which the symbol $(\infty)$ 
indicates that there are infinitely many edges from $v$ to $w$.   

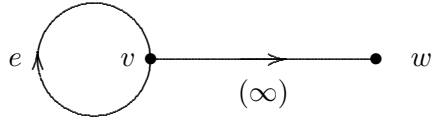
\[ \beginpicture
\setcoordinatesystem units <1.5cm,1.5cm>
\setplotarea x from -7 to 10, y from -0.7 to 0.5
\put {$\bullet$} at -4 0
\put {$\bullet$} at -2 0

\put {$v$} [v] at -4.2 0 
\put {$w$} [v] at -1.6 0 
\put {$e$}     at -5.2 0

\put {$(\infty)$} at -3 -0.3 

\setlinear 
\plot -4 0  -2 0 /  

\arrow <0.235cm> [0.2,0.6] from -5 0 to -5 0.07 
\arrow <0.235cm> [0.2,0.6] from -3.2 0 to -2.8 0  

\circulararc 360 degrees from -4 0 center at -4.5 0

\endpicture \] 
There are two maximal tails in $E_1$, namely $E_1^0$ and 
$M=\{v\}$. $E_1^0$ is in ${\mathcal M}_\gamma(E_1)$ and 
$M$  belongs to ${\mathcal M}_\tau(E_1)$. There is a 
unique breaking vertex $v$ in $E_1$. The bijection of 
Corollary \ref{primitiveideals} identifies $\{E_1^0,\,v\}
\cup(M\times\T)$ with $\Prim(C^*(E_1))$. The topology can 
be determined by Theorem \ref{9closures}. 
The closure of $\{E_1^0\}$ is the 
entire space $\Prim(C^*(E_1))$, the closure of $\{v\}$ is 
$\{v\}\cup(M\times\T)$, and for every $D\subseteq\T$ the 
closure of $M\times D$ is $M\times\overline{D}$. 

The maximal tail $E_1^0$ corresponds to the primitive ideal 
$\{0\}$. The breaking vertex $v$ corresponds to 
the ideal $I_w$, generated by the projection $p_w$. This ideal  
is isomorphic with the compacts 
and is essential in $C^*(E_1)$ by \cite[Lemma 1.1]{sz2}. The 
quotient $C^*(E_1)/I_w$ is isomorphic to the Toeplitz 
algebra ${\mathcal T}$ by \cite[Proposition 3.4]{bhrs}. 
Thus, there is a short exact sequence 
$$ 0\longrightarrow\K\longrightarrow C^*(E_1)\longrightarrow 
   {\mathcal T}\longrightarrow 0. $$ 
As shown in Lemma \ref{sandwich}, each non gauge-invariant 
primitive ideal $R_{M,t}$, corresponding to the maximal tail 
$M$ and $t\in\T$, is sandwiched between two gauge-invariant ideals, 
namely 
$$ J_{\{w\},\{v\}}\subset R_{M,t}\subset C^*(E_1). $$ 
The ideal $J_{\{w\},\{v\}}$ is generated by the projections 
$p_w$ and $p_v-s_es_e^*$. 
\end{example} 

\begin{example}\label{example2}
The following graph $E_2$, considered in \cite{fmr}, is neither 
row-finite nor does it satisfy Condition (K). 

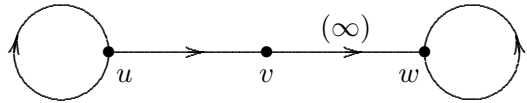
\[ \beginpicture
\setcoordinatesystem units <2.1cm,2.1cm>
\setplotarea x from -3 to 2, y from -0.4 to 0.1
\put {$\bullet$} at -2 0
\put {$\bullet$} at -1 0
\put {$\bullet$} at 0 0
\put {$(\infty)$} at -0.5 0.15

\setlinear
\plot -2 0 0 0 /
\arrow <0.235cm> [0.2,0.6] from -1.5 0 to -1.4 0
\arrow <0.235cm> [0.2,0.6] from -0.5 0 to -0.4 0

\circulararc 360 degrees from -2 0 center at -2.3 0
\circulararc 360 degrees from 0 0 center at 0.3 0

\arrow <0.235cm> [0.2,0.6] from -2.6 0 to -2.58 0.1
\arrow <0.235cm> [0.2,0.6] from 0.6 0 to 0.58 0.1

\put {$u$} [v] at -1.9  -0.15 
\put {$v$} [v] at -1 -0.15
\put {$w$} [v] at -0.1 -0.15

\endpicture \]
There are no breaking vertices in $E_2$. There are three maximal tails in $E_2$: 
$M_1=\{u\}$, $M_2=\{u,v\}$ and $M_3=\{u,v,w\}$. $M_2\in{\mathcal M}_\gamma(E_2)$, 
while $M_1$ and $M_3$ belong to ${\mathcal M}_\tau(E_2)$. The bijection 
of Corollary \ref{primitiveideals} identifies $\{M_2\}\cup(M_1\times\T) 
\cup(M_3\times\T)$ with $\Prim(C^*(E_2))$. The topology is given 
by Theorem \ref{9closures}. The closure of $\{M_2\}$ is 
$\{M_2\}\cup(M_1\times\T)$. For any $D\subseteq\T$, the closure of 
$M_1\times D$ is $M_1\times\overline{D}$, and the closure of 
$M_3\times D$ is $\{M_2\}\cup(M_1\times\T)\cup(M_3\times\overline{D})$. 
As in Lemma \ref{sandwich}, for any $t\in\T$ we have 
$$ I_{\{v,w\}}\subset R_{M_1,t}\subset C^*(E_2)\;\; \text{ and } \;\; 
   \{0\}\subset R_{M_3,t}\subset I_{w}. $$ 
\end{example}

\begin{example}\label{example3}
Let $E_3$ be the following graph with $E_3^0=\{v_{i, j}\, |\, 1\leq i, j <
\infty\}$ and $E_3^1=\{e_i\}\cup\{f_{i, j}\}\cup\{g_{i,j}\}$, where
\begin{eqnarray*}
s(e_i)=v_{1, 2i}, & s(f_{i,j})=v_{i,j}, & s(g_{i,j})=v_{i,j}, \\
r(e_i)=v_{1, 2i-1}, & r(f_{i,j})=v_{i,j+1}, & r(g_{i,j})=v_{i+1,j}.
\end{eqnarray*}

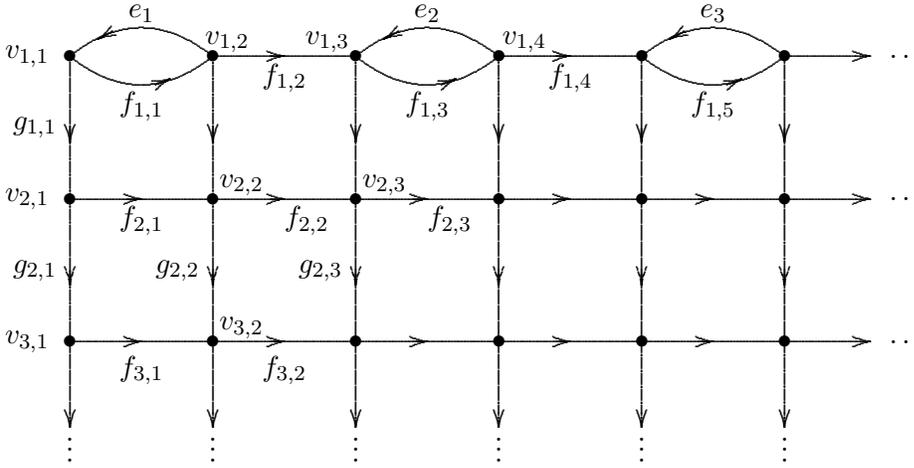
\[ \beginpicture
\setcoordinatesystem units <1.9cm,1.9cm>
\setplotarea x from -3 to 5, y from -1 to 2 
\put {$\bullet$} at -3 0
\put {$\bullet$} at -2 0
\put {$\bullet$} at -1 0
\put {$\bullet$} at 0 0
\put {$\bullet$} at 1 0
\put {$\bullet$} at 2 0

\put {$\bullet$} at -3 1
\put {$\bullet$} at -2 1
\put {$\bullet$} at -1 1
\put {$\bullet$} at 0 1
\put {$\bullet$} at 1 1
\put {$\bullet$} at 2 1

\put {$\bullet$} at -3 2
\put {$\bullet$} at -2 2
\put {$\bullet$} at -1 2
\put {$\bullet$} at 0 2
\put {$\bullet$} at 1 2
\put {$\bullet$} at 2 2

\setquadratic
\plot -3 2 -2.5 2.2 -2 2 /
\plot -1 2 -0.5 2.2 0 2 /
\plot 1 2 1.5 2.2 2 2 /
\plot -3 2 -2.5 1.8 -2 2 /
\plot -1 2 -0.5 1.8 0 2 /
\plot 1 2 1.5 1.8 2 2 /

\setlinear
\plot -3 2 -3 -0.5 /
\plot -2 2 -2 -0.5 /
\plot -1 2 -1 -0.5 /
\plot 0 2 0 -0.5  /
\plot -2 2 -1 2 /
\plot 0 2 1 2 /
\plot 2 2 2.5 2 /
\plot -3 1 2.5 1 /
\plot -3 0 2.5 0 /
\plot 1 2 1 -0.5 /
\plot 2 2 2 -0.5 /

\put {$v_{1,1}$} [v] at -3.3 2
\put {$v_{2,1}$} [v] at -3.3 1
\put {$v_{3,1}$} [v] at -3.3 0

\put {$v_{1,2}$} [v] at -1.9 2.1
\put {$v_{1,3}$} [v] at -1.2 2.1
\put {$v_{1,4}$} [v] at 0.18 2.1

\put {$v_{2,2}$} [v] at -1.8 1.1
\put {$v_{3,2}$} [v] at -1.8 0.1
\put {$v_{2,3}$} [v] at -0.8 1.1

\put {$g_{1,1}$} [l] at -3.4 1.5
\put {$g_{2,1}$} [l] at -3.4 0.5
\put {$g_{2,2}$} [l] at -2.4 0.5
\put {$g_{2,3}$} [l] at -1.4 0.5

\put {$e_1$} [v] at -2.5 2.3
\put {$e_2$} [v] at -0.5 2.3
\put {$e_3$} [v] at 1.5 2.3

\put {$f_{1,2}$} [v] at -1.5 1.85
\put {$f_{1,4}$} [v] at 0.5 1.85
\put {$f_{1,1}$} [v] at -2.5 1.65
\put {$f_{1,3}$} [v] at -0.5 1.65
\put {$f_{1,5}$} [v] at 1.5 1.65

\put {$f_{2,1}$} [v] at -2.5 0.85
\put {$f_{3,1}$} [v] at -2.5 -0.2
\put {$f_{3,2}$} [v] at -1.5 -0.2

\put {$f_{2,2}$} [l] at -1.5 0.85
\put {$f_{2,3}$} [l] at -0.5 0.85

\arrow <0.235cm> [0.2,0.6] from -2.6 1 to -2.5 1
\arrow <0.235cm> [0.2,0.6] from -1.6 1 to -1.5 1
\arrow <0.235cm> [0.2,0.6] from -0.6 1 to -0.5 1
\arrow <0.235cm> [0.2,0.6] from 0.4 1 to 0.5 1
\arrow <0.235cm> [0.2,0.6] from 1.4 1 to 1.5 1

\arrow <0.235cm> [0.2,0.6] from -2.6 0 to -2.5 0
\arrow <0.235cm> [0.2,0.6] from -1.6 0 to -1.5 0
\arrow <0.235cm> [0.2,0.6] from -0.6 0 to -0.5 0
\arrow <0.235cm> [0.2,0.6] from 0.4 0 to 0.5 0
\arrow <0.235cm> [0.2,0.6] from 1.4 0 to 1.5 0

\arrow <0.235cm> [0.2,0.6] from -3 1.5 to -3 1.4
\arrow <0.235cm> [0.2,0.6] from -3 0.5 to -3 0.4
\arrow <0.235cm> [0.2,0.6] from -3 -0.5 to -3 -0.6

\arrow <0.235cm> [0.2,0.6] from -2 1.5 to -2 1.4
\arrow <0.235cm> [0.2,0.6] from -2 0.5 to -2 0.4
\arrow <0.235cm> [0.2,0.6] from -2 -0.5 to -2 -0.6

\arrow <0.235cm> [0.2,0.6] from -1 1.5 to -1 1.4
\arrow <0.235cm> [0.2,0.6] from -1 0.5 to -1 0.4
\arrow <0.235cm> [0.2,0.6] from -1 -0.5 to -1 -0.6

\arrow <0.235cm> [0.2,0.6] from 0 1.5 to 0 1.4
\arrow <0.235cm> [0.2,0.6] from 0 0.5 to 0 0.4
\arrow <0.235cm> [0.2,0.6] from 0 -0.5 to 0 -0.6

\arrow <0.235cm> [0.2,0.6] from 1 1.5 to 1 1.4
\arrow <0.235cm> [0.2,0.6] from 1 0.5 to 1 0.4
\arrow <0.235cm> [0.2,0.6] from 1 -0.5 to 1 -0.6

\arrow <0.235cm> [0.2,0.6] from 2 1.5 to 2 1.4
\arrow <0.235cm> [0.2,0.6] from 2 0.5 to 2 0.4
\arrow <0.235cm> [0.2,0.6] from 2 -0.5 to 2 -0.6

\arrow <0.235cm> [0.2,0.6] from 2.5 2 to 2.6 2
\arrow <0.235cm> [0.2,0.6] from 2.5 1 to 2.6 1
\arrow <0.235cm> [0.2,0.6] from 2.5 0 to 2.6 0

\arrow <0.235cm> [0.2,0.6] from -1.6 2 to -1.5 2
\arrow <0.235cm> [0.2,0.6] from 0.4 2 to 0.5 2

\arrow <0.235cm> [0.2,0.6] from -2.7 2.175 to -2.8 2.12
\arrow <0.235cm> [0.2,0.6] from -2.4 1.8 to -2.3 1.84

\arrow <0.235cm> [0.2,0.6] from -0.7 2.175 to -0.8 2.12
\arrow <0.235cm> [0.2,0.6] from -0.4 1.8 to -0.3 1.84

\arrow <0.235cm> [0.2,0.6] from 1.3 2.175 to 1.2 2.12
\arrow <0.235cm> [0.2,0.6] from 1.6 1.8 to 1.7 1.84

\put {$\vdots$} at -3 -0.7
\put {$\vdots$} at -2 -0.7
\put {$\vdots$} at -1 -0.7
\put {$\vdots$} at 0 -0.7
\put {$\vdots$} at 1 -0.7
\put {$\vdots$} at 2 -0.7

\put {$\ldots$} at 2.85 2
\put {$\ldots$} at 2.85 1
\put {$\ldots$} at 2.85 0
\endpicture \]
This is an infinite row-finite graph which does not satisfy
Condition (K). Since $E_3$ is row-finite there are 
no breaking vertices. There are four families of maximal tails, 
indexed by the integers $n\geq 1$:
\begin{align*}
M_n & =  \{v_{i,j}:1\leq i\leq n,\:1\leq j<\infty\}, \\
M^{2n-1} & =  \{v_{i,j}:1\leq i<\infty,\:
   1\leq j\leq 2n-1\}\cup\{v_{1,2n}\}, \\
M^{2n} & =  \{v_{i,j}:1\leq i<\infty,\:1\leq j\leq 2n\}, \\
T_n & =  \{v_{1,j}:1\leq j\leq 2n\}.
\end{align*}
In addition, $E_3^0$ is a maximal tail too. $E_3^0$ and all $M_n$ and $M^n$ 
belong to ${\mathcal M}_\gamma (E_3)$. On the other hand, 
each maximal tail $T_n$ contains a loop without exits and 
hence $T_n\in{\mathcal M}_\tau (E_3)$. By Corollary \ref{primitiveideals}, 
there is a bijection between $\{E_3^0\}\cup
\{M_n:n\geq1\}\cup\{M^n:n\geq1\}\cup\bigcup_{n\geq1}(T_n\times\T)$ 
and $\Prim(C^*(E_3))$. The topology of $\Prim(C^*(E_3))$ can 
be determined with help of Corollary \ref{4closures}. 
\end{example}

\end{document}